\documentclass[12pt]{article}

\usepackage{amsmath}
\usepackage{amsfonts}
\usepackage{graphicx}
\usepackage{rotating}

\RequirePackage[numbers]{natbib}

\renewcommand{\span}{{\rm span}}
\newcommand{\inverse}[1]{{\textstyle\frac{1}{#1}}}
\newcommand{\half}{\inverse{2}}

\newcommand{\R}{{\mathbb R}}
\newcommand{\CPP}{{\ttfamily C++ }}
\newcommand{\dif}{{\mathrm d}}
\DeclareMathOperator{\logit}{logit}
\newcommand{\class}[1]{{\sffamily{#1}}}
\newcommand{\method}[1]{{\sffamily{#1}}}

\title{\vspace{-4cm} {\Large \bf Stochastic filtering
via $L^2$ projection on mixture manifolds
with computer algorithms and numerical examples}\\ {\small Updated version in Mathematics of Control, Signals \& Systems 28(1), 1-33, 2016}}

\author{John Armstrong  \\
Dept. of Mathematics, \ 
King's College, London\\
{\tt john.1.armstrong@kcl.ac.uk}\\
\\
Damiano Brigo  \\
Dept. of Mathematics, \ \ 
Imperial College, London\\
{\tt damiano.brigo@imperial.ac.uk}}

\date{\small {First version: March 26, 2013. This version: \today }}
\newtheorem{theorem}{Theorem}[section]

\newtheorem{example}[theorem]{Example}

\begin{document}

\maketitle
\thispagestyle{empty}

\begin{abstract}
We examine some differential geometric approaches to
finding approximate solutions to the continuous time nonlinear
filtering problem. Our primary focus is a new projection method for the optimal filter infinite dimensional Stochastic Partial Differential Equation (SPDE), based on the
direct L2 metric and on a family of normal mixtures. We compare this method to 
earlier projection methods based on the Hellinger distance/Fisher metric and exponential families, 
and we compare the L2 mixture projection filter with a particle method with the same number of parameters, using the Levy metric. We prove that for a simple choice of the mixture manifold the L2 mixture projection filter coincides with a Galerkin method, whereas for more general mixture manifolds the equivalence does not hold and the L2 mixture filter is more general. 
We study particular systems that may illustrate the advantages of this new filter over other algorithms  when comparing outputs with the optimal filter. 
We finally consider a specific software design that is suited for a numerically efficient implementation of this filter and provide numerical examples. 
\end{abstract}

\bigskip


%\medskip

\vspace{0.2cm}

\noindent { {\bf AMS} Classification codes: 53B25, 53B50, 60G35, 62E17, 62M20, 93E11}

\pagestyle{myheadings}
\markboth{} {J. Armstrong, D. Brigo. Stochastic
Filtering via $L^2$ Projection}

\tableofcontents

\newpage

\section{Introduction}

In the nonlinear filtering problem one observes a system whose state
is known to follow a given stochastic differential equation. The
observations that have been made contain an additional noise term, so
one cannot hope to know the true state of the system. However, one can
reasonably ask what is the probability density over the possible states.

When the observations are made in continuous time, the probability
density follows a stochastic partial differential equation known as
the Kushner--Stratonovich equation. This can be seen as a generalization
of the Fokker--Planck equation that expresses the evolution of the
density of a diffusion process. Thus the problem we wish to address boils
down to finding approximate solutions to the Kushner--Stratonovich equation.

For a quick introduction to the filtering problem see Davis and Marcus
(1981) \citep{davis81b}. For a more complete treatment from a
mathematical point of view see Lipster and Shiryayev (1978)
\citep{LipShi} and Bain and Crisan~\citep{crisan10}. See Jazwinski (1970) \citep{jazwinski70a} and Ahmed (1998) \cite{ahmedbook} for
a comprehensive treatment of filtering. 
For recent results see the volume  \citep{crisan}.

The main idea we will employ is inspired by the differential geometric
approach to statistics developed in \citep{amari85a} and
\citep{pistonesempi}. The idea of applying this approach to the filtering problem has been sketched first in \citep{hanzon87}. 
One thinks of the probability distribution as evolving in an infinite dimensional space ${\cal P}$
which is in turn contained in some Hilbert space $H$.  One can then
think of the Kushner--Stratonovich equation as defining a vector field
in ${\cal P}$: the integral curves of the vector field should
correspond to the solutions of the equation.  To find approximate
solutions to the Kushner--Stratonovich equation one chooses a finite
dimensional submanifold $M$ of $H$ and approximates the probability
distributions as points in $M$. At each point of $M$ one can use the
Hilbert space structure to project the vector field onto the tangent
space of $M$. One can now attempt to find approximate solutions to the
Kushner--Stratonovich equations by integrating this vector field on
the manifold $M$.

This mental image is slightly innaccurate. The Kushner--Stratonovich
equation is a stochastic PDE rather than a PDE so one should imagine
some kind of stochastic vector field rather than a smooth vector
field. Thus in this approach we hope to approximate the infinite
dimensional stochastic PDE by solving a finite dimensional stochastic
ODE on the manifold.

Note that our approximation will depend upon two choices: the choice
of manifold $M$ and the choice of Hilbert space structure $H$. In this
paper we will consider two possible choices for the Hilbert space
structure: the direct $L^2$ metric on the space of probability
distributions; the Hilbert space structure associated with the
Hellinger distance and the Fisher Information metric. Our focus will
be on the direct $L^2$ metric since projection using the Hellinger
distance has been considered before. As we shall see, the choice of
the ``best'' Hilbert space structure is determined by the manifold one
wishes to consider --- for manifolds associated with exponential
families of distributions the Hellinger metric leads to the simplest
equations, whereas the direct $L^2$ metric works well with mixture
distributions. 

It was proven in \cite{brigo99} that the projection filter in Hellinger metric on exponential families is equivalent to the classical assumed density filters. In this paper, we show that the projection filter for basic mixture manifolds in L2 metric is equivalent to a Galerkin method. This only holds for very basic mixture families, however, and the L2 projection method turns out to be more general than the Galerkin method. 

We will write down the stochastic ODE determined by the geometric
approach when $H=L^2$ and show how it leads to a numerical scheme for
finding approximate solutions to the Kushner--Stratonovich equations
in terms of a mixture of normal distributions. We will call this
scheme the {\em $L^2$ normal mixture projection filter} or simply
the L2NM projection filter. 

The stochastic ODE for the Hellinger metric was considered in
\citep{BrigoPhD}, \citep{brigo98} and
\citep{brigo99}.  In particular a precise numerical scheme
is given in \citep{BrigoPhD} for finding solutions by projecting onto
an exponential family of distributions. We will call this scheme the
{\em Hellinger exponential projection filter} or simply the HE projection
filter.

En passant, we clarify why we use the Kushner--Stratonovich and not the alternative Duncan--Mortensen--Zakai (Zakai for brevity) equation for nonlinear filtering. The Zakai equation is the preferred stochastic 
partial differential equation in nonlinear filtering \cite{ahmedbook}, since it is a linear equation and admits a robust version. It models a particular unnormalized version of the optimal filter density. We will show that projecting the Zakai equation in Hellinger metric results in the same approximate filter as projecting the Zakai equation. Intuitively, in Hellinger metric square roots of densities are part of the unit sphere in $L^2$ and re-scaling is orthogonal to the tangent space. Thus in Hellinger metric we could indeed use the Zakai equation and we would obtain the same filter. However, when projecting according to the direct metric, this no longer holds and projecting the Zakai equation would lead to a different filter than projecting the Kushner Stratonovich equation. Since we are interested in approximating well the density and not an unnormalized version of it, we use the Kushner--Stratonovich equation throughout the paper.

We will compare the results of a \CPP
implementation of the L2NM projection filter with a number
of other numerical approaches including the HE projection
filter and the optimal filter. We can measure the goodness of our filtering approximations 
thanks to the geometric structure and, in particular, the precise metrics we are using 
on the spaces of probability measures. 

What emerges is that the two projection methods produce excellent results for
a variety of filtering problems. The results appear similar for both
projection methods; which gives more accurate results depends upon the
problem.

As we shall see, however, the L2NM projection approach can
be implemented more efficiently. In particular one needs to perform
numerical integration as part of the HE projection filter algorithm
whereas all integrals that occur in the L2NM projection can be
evaluated analytically.

We also compare the L2NM filter to a particle filter with the best possible combination of particles
with respect to the L{\'e}vy metric. Introducing the L{\'e}vy metric is needed because particles densities do not compare well with smooth densities when using $L^2$ induced metrics.
We show that, given the same number of parameters, the L2NM may outperform a particles based system. 

The paper is structured as follows: In Section \ref{sec:nonlinfp} we introduce the nonlinear filtering problem and the infinite-dimensional Stochastic PDE (SPDE) that solves it.  In Section \ref{sec:SMan} we introduce the geometric structure we need to project the filtering SPDE onto a finite dimensional manifold of probability densities. In Section \ref{sec:PF} we perform the projection of the filtering SPDE according to the L2NM framework and also recall the HE based framework. 
In Section \ref{sec:Galerk} we prove equivalence between the projection filter in L2 metric for basic mixture families and the Galerkin method. 
In Section \ref{sec:NImp} we briefly introduce the main issues in the numerical implementation and then focus on software design for the L2NM filter. 
In Section \ref{sec:normmixcase} a second theoretical result is provided, showing a particularly convenient structure for the projection filter equations for specific choices of the system properties and of the mixture manifold. 

In Section \ref{NumRes} we look at numerical results, whereas in Section \ref{sec:Part} we compare our outputs with a particle method. Section \ref{sec:Conc} concludes the paper.

\section{The non-linear filtering problem with\\ continuous-time  observations}\label{sec:nonlinfp}

In the non-linear filtering problem the state of some system is
modelled by a process $X$ called the signal. This signal evolves over
time $t$ according to an It\^{o} stochastic differential equation
(SDE). We measure the state of the system using some observation
$Y$. The observations are not accurate, there is a noise term. So the
observation $Y$ is related to the signal $X$ by a second equation.

\begin{equation} \label{Lanc1-1}
\begin{array}{rcl}
   dX_t &=& f_t(X_t)\,dt + \sigma_t(X_t)\,dW_t, \ \ X_0, \\ \\
   dY_t &=& b_t(X_t)\,dt + dV_t, \ \ Y_0 = 0\ .
\end{array}
\end{equation}

In these equations the unobserved state process $\{ X_t, t \geq 0 \}$
takes values in $\R^n$, the observation $\{ Y_t, t\geq 0 \}$ takes
values in $\R^d$ and the noise processes $\{ W_t, t\geq 0\}$ and $\{
V_t, t\geq 0\}$ are two Brownian motions.

The nonlinear filtering problem consists in finding the conditional
probability distribution $\pi_t$ of the state $X_t$ given the
observations up to time $t$ and the prior distribution $\pi_0$ for $X_0$.

Let us assume that $X_0$, and the two Brownian motions are
independent. Let us also assume that the covariance matrix for $V_t$
is invertible. We can then assume without any further loss of
generality that its covariance matrix is the identity. We
introduce a variable $a_t$ defined by:

\[ a_t = \sigma_t \sigma_t^T \]

With these preliminaries, and a number of rather more technical
conditions which we will state shortly, one can show that $\pi_t$
satisfies the a stochastic PDE called the Kushner--Stratonovich
equation. This states that for any compactly supported test function
$\phi$ defined on $\R^{n}$

\begin{equation} \label{FKK}
   \pi_t(\phi) = \pi_0(\phi)
   + \int_0^t \pi_s({\cal L}_s \phi)\,ds
   + \sum_{k=1}^d \int_0^t
   [\pi_s(b_s^k\, \phi) - \pi_s(b_s^k)\, \pi_s(\phi)]\,
   [dY_s^k-\pi_s(b_s^k)\,ds]\ ,
\end{equation}
where for all $t \geq 0$,
the backward diffusion operator ${\cal L}_t$ is defined by
\begin{displaymath}
   {\cal L}_t = \sum_{i=1}^n f_t^i\,
   \frac{\partial}{\partial x_i} + \half \sum_{i,j=1}^n
   a_t^{ij}\, \frac{\partial^2}{\partial x_i \partial x_j}\ .
\end{displaymath}

Equation~(\ref{FKK}) involves the derivatives of the test function
${\phi}$ because of the expression ${\cal L}_s \phi$. We assume now
that $\pi_t$ can be represented by a density $p_t$ with respect to the
Lebesgue measure on $\R^n$ for all time $t \geq 0$ and that
we can replace the term involving ${\cal L}_s\phi$  with a term involving its
formal adjoint ${\cal L}^*$. Thus, proceeding formally, we find that $p_t$ obeys
the following It\^o-type stochastic partial differential equation (SPDE):
\[ \dif p_t = {\cal L}^*_t p_t \dif t + \sum_{k=1}^d p_t [ b_t^k -
E_{p_t} \{ b_t^k \} ][ \dif Y_t^k - E_{p_t} \{b_t^k \} \dif t ] \]
where $E_{p_t}\{ \cdot \}$ denotes the expectation with respect to the
probability density $p_t$ (equivalently the conditional expectation
given the observations up to time $t$). The forward diffusion operator ${\cal L}^*_t$
is defined by:
\begin{displaymath}
   {\cal L}_t^* \phi = -\sum_{i=1}^n 
   \frac{\partial}{\partial x_i} [ f_t^i \phi ] + \half \sum_{i,j=1}^n
    \frac{\partial^2}{\partial x_i \partial x_j}[ a_t^{ij} \phi ].
\end{displaymath}

This equation is written in It\^{o} form. When working with stochastic
calculus on manifolds it is necessary to use Stratonovich SDE's rather
than It\^{o} SDE's. This is because one does not in general know how
to interpret the second order terms that arise in It\^{o} calculus in
terms of manifolds. The interested reader should consult
\citep{elworthy82a}. A straightforward calculation yields the following
Stratonvich SPDE:

\begin{displaymath}
   dp_t = {\cal L}_t^\ast\, p_t\,dt
   - \half\, p_t\, [\vert b_t \vert^2 - E_{p_t}\{\vert b_t \vert^2\}] \,dt
   + \sum_{k=1}^d p_t\, [b_t^k-E_{p_t}\{b_t^k\}] \circ dY_t^k\ .
\end{displaymath}

We have indicated that this is the Stratonovich form of the equation
by the presence of the
symbol `$\circ$' inbetween the diffusion coefficient and the Brownian
motion of the SDE. We shall use this convention throughout the rest of
the paper.

In order to simplify notation, we introduce the following
definitions~:
\begin{equation} \label{coeff}
\begin{array}{rcl}
 \gamma_t^0(p) &:=&
   \half\, [\vert b_t \vert^2 - E_p\{\vert b_t \vert^2\}]\ p, \\ \\
   \gamma_t^k(p) &:=& [b_t^k - E_p\{b_t^k\}] p \ ,
\end{array}
\end{equation}
for $k = 1,\cdots,d$.
The Stratonovich form of the Kushner--Stratonovich equation reads now
\begin{equation}\label{KSE:str}
   dp_t = {\cal L}_t^\ast\, p_t\,dt
   -  \gamma_t^0(p_t)\,dt
   + \sum_{k=1}^d  \gamma_t^k(p_t) \circ dY_t^k\ .
\end{equation}

Thus, subject to the assumption that a density $p_t$ exists for all
time and assuming the necessary decay condition to ensure that
replacing ${\cal L}$ with its formal adjoint is valid, we find that
solving the non-linear filtering problem is equivalent to solving this
SPDE. Numerically approximating the solution of equation~(\ref{KSE:str})
is the primary focus of this paper.

For completeness we review the technical conditions required in order
for equation~(\ref{FKK}) to follow from~(\ref{Lanc1-1}). 

\begin{itemize}
   \item[(A)] Local Lipschitz continuity~:
for all $R > 0$, there exists $K_R > 0$ such that
\begin{displaymath}
   \vert f_t(x) - f_t(x') \vert \leq K_R\, \vert x-x' \vert
   \hspace{1cm}\mbox{and}\hspace{1cm}
   \Vert a_t(x) - a_t(x') \Vert \leq K_R\, \vert x-x' \vert\ ,
\end{displaymath}
for all $t\geq 0$, and for all $x,x'\in B_R$, the ball of radius $R$.

   \item[(B)] Non--explosion~:
there exists $K > 0$ such that
\begin{displaymath}
   x^T f_t(x) \leq K\, (1+\vert x \vert^2)
   \hspace{1cm}\mbox{and}\hspace{1cm}
   \mbox{trace}\; a_t(x) \leq K\, (1+\vert x \vert^2)\ ,
\end{displaymath}
for all $t\geq 0$, and for all $x\in {\bf R}^n$.

   \item[(C)] Polynomial growth~:
there exist $K > 0$ and $r \geq 0$ such that
\begin{displaymath}
   \vert b_t(x) \vert \leq K\, (1+\vert x \vert^r)\ ,
\end{displaymath}
for all $t\geq 0$, and for all $x\in {\bf R}^n$.
\end{itemize}

Under assumptions~(A) and~(B), there exists
a unique solution $\{X_t\,,\,t\geq 0\}$ to the state equation,
see for example~\citep{khasminskii},
and $X_t$ has finite moments of any order.
Under the additional assumption~(C) the following {\em finite energy}
condition holds
\begin{displaymath}
   E \int_0^T \vert b_t(X_t) \vert^2\,dt < \infty\ ,
   \hspace{1cm}\mbox{for all $T\geq 0$}.
\end{displaymath}

Since the finite energy condition holds, it follows from Fujisaki,
Kallianpur and Kunita~\citep{fujisaki72a} that $\{\pi_t\,,\,t\geq 0\}$
satisfies the Kushner--Stratonovich equation~(\ref{FKK}).

In closing this summary of nonlinear filtering, we should point out that, in a broad part of the nonlinear filtering literature, the preferred SPDE for the optimal filter is a SPDE for an unnormalized version $q$ of the optimal filter density $p$. The Duncan-Mortensen-Zakai equation (Zakai for brevity) for the unnormalized density $q_t(x)$ of the optimal filter
%, such that
%\[ p_t(x) = \frac{q_t(x)}{ \int q_t(x) dx} \]
reads, in Stratonovich form 
 \begin{displaymath}
   dq_t = {\cal L}_t^\ast\, q_t\,dt
   - \half\, q_t\, \vert b_t \vert^2  \ dt
   + \sum_{k=1}^d q_t\, [b_t^k] \circ dY_t^k\ , \ \ q_0 = p_0,
\end{displaymath}
see for example Eq. 14.31 in \cite{ahmedbook}, where conditions under which this is an evolution equation in $L^2$ are discussed.
This is a linear Stochastic PDE and as such it is more tractable than the Kushner-Stratonovich (KS) Equation. Linearity has led the Zakai Eq. to be preferred to the KS Eq. in general. The reason why we still resort to KS will be clarified in full detail when we derive the projection filter below, although we may anticipate that this has to do with the fact that we aim to derive an approximation that is good for the normalized density $p$ and not for the unnormalized density $q$. 
There is a second and related reason why the Zakai Eq. has been preferred to the KS Eq. in the literature: the possibility to derive a robust PDE for the optimal filter, which we briefly review now.

While here we use Stratonovich calculus in order to deal with manifolds, nonlinear filtering equations are usually based on Ito stochastic differential equations, holding almost surely in the space of continuous functions. However, in practice real-world sample paths have finite variation, and the set of finite variation paths has measure zero under the Wiener measure. It is in principle possible to obtain a version of the filter which takes arbitrary values on any real-world sample path. In technical terms, nonlinear filtering equations such as the KS or Zakai equation are not robust. One would wish to have that the nonlinear filter is a continuous functional on the continuous functions (paths). This way we would have that the filter based on real-life finite variation paths is close to the theoretical optimal one based on unbounded variation paths. There are essentially two approaches to bypass the above lack of robustness. In one approach, introduced by Clark \cite{clark78}, the Zakai stochastic PDE is transformed into an equivalent pathwise form avoiding the differential $d Y_t$ of the observation process. Indeed, assuming the observation function $b$ to be time homogeneous for simplicity, $b_t(x) = b(x)$, set 
\[ r_t(x) :=  \exp\left(-\sum_{k=1}^d b^k(x) Y^k_t \right) q_t(x) . \] 
Then by the standard chain rule of Stratonovich calculus one derives easily 
\[ \partial_t  r_t =  \frac{{\cal L}_t^\ast\ \left( \exp\left(\sum_{k=1}^d b^k Y^k_t \right) r_t \right)}{\exp\left(\sum_{k=1}^d b^k Y^k_t \right) }
   - \half\,  \vert b \vert^2\ r_t .  
\]
This transformed filtering equation is shown to be continuous  with respect to the observation under a suitable topology. This equation can then be extended  to real-life sample paths using continuity \cite{davis80}. 

In the other approach to robust filtering, due originally to Balakrishnan \cite{balakrishnan80}, one tries to model the observation process directly with a white noise error term. Although modelling white noise directly is intuitively appealing, it brings a host of mathematical complications. Kallianpur and Karandikar \cite{kallianpur83} developed the theory of nonlinear filtering in this framework, while \cite{bagchi94} extend this second approach to the case of correlated state and observation noises. 

Going back to the projection filter, one could consider using a projection method to find approximate solutions for $r_t$. However, as our primary interest is in finding good approximations to $p_t$ rather than $r_t$, we believe that projecting the KS equation is the more promising approach.

\section{Statistical manifolds }\label{sec:SMan}  
\subsection{ Families of distributions }

As discussed in the introduction, the idea of a projection filter is
to approximate solutions to the Kushner--Stratononvich equation~(\ref{FKK})
using a finite dimensional family of distributions.

\begin{example}
A {\em normal mixture} family contains distributions given by:
\[ p = \sum_{i=1}^m \lambda_i \frac{1}{\sigma_i \sqrt{ 2 \pi}}
\exp\left(\frac{-(x-\mu_i)^2}{2 \sigma_i^2}\right) \]
with $\lambda_i >0 $ and $\sum \lambda_i = 1$. It is
a $3m-1$ dimensional family of distributions.
\end{example}

\begin{example}
A {\em polynomial exponential family} contains distributions given by:
\[ p = \exp( \sum_{i=0}^m a_i x^i ) \] where $a_0$ is chosen to ensure
that the integral of $p$ is equal to $1$.  To ensure the convergence
of the integral we must have that $m$ is even and $a_m < 0$. This is
an $m$ dimensional family of distributions. Polynomial exponential families
are a special case of the more general notion of an exponential family, see
for example \citep{amari85a}.
\end{example}

A key motivation for considering these families is that one can
reproduce many of the qualitative features of distributions that arise
in practice using these distributions. For example, consider the
qualitative specification: the distribution should be bimodal with
peaks near $-1$ and $1$ with the peak at $-1$ twice as high and twice
as wide as the peak near $1$. One can easily write down a distribution of
this approximates form using a normal mixture.

\begin{figure}[htp]
\begin{centering}
\includegraphics{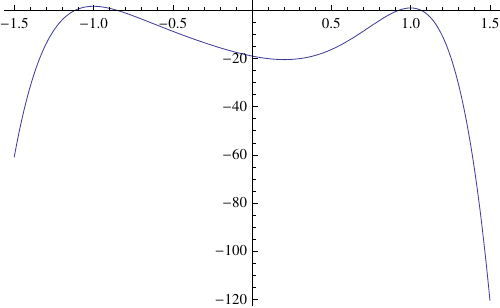}
\caption{ $y = -18.98-13.15 x+23.54 x^2+25.43 x^3+13.96 x^4-12.63 x^5-17.15 x^6$}
\label{fig:bimodalsextic}
\end{centering}
\end{figure}

\begin{figure}[htp]
\begin{centering}
\includegraphics{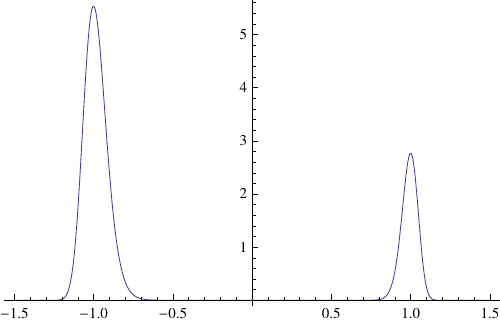}
\caption{ $y = \exp(-18.98-13.15 x+23.54 x^2+25.43 x^3+13.96 x^4-12.63 x^5-17.15 x^6)$}
\label{fig:bimodalexponential}
\end{centering}
\end{figure}

To find a similar exponential family, one seeks a polynomial with:
local maxima at $-1$ and $1$; with the maximum values at these points
differing by $\log(2)$; with second derivative at $1$ equal to twice
that at $-1$. These conditions give linear equations in the polynomial
coefficients. Using degree $6$ polynomials it is simple to find
solutions meeting all these requirements. A specific numerical example
of a polynomial meeting these requirements is plotted in
Figure~\ref{fig:bimodalsextic}. The associated exponential
distribution is plotted in Figure~\ref{fig:bimodalexponential}.

We see that normal mixtures and exponential families have a broadly
similar power to describe the qualitative shape of a distribution
using only a small number of parameters. Our hope is that by
approximating the probability distributions that occur in the
Kushner--Stratonovich equation by elements of one of these families we
will be able to derive a low dimensional approximation to the full
infinite dimensional stochastic partial differential equation.

\subsection{Two Hilbert spaces of probability distributions}

We have given direct parameterisations of our families of probability
distributions and thus we have implicitly represented them as finite
dimensional manifolds. In this section we will see how families of
probability distributions can be thought of as being embedded in a
Hilbert space and hence they inherit a manifold structure and metric
from this Hilbert space.

There are two obvious ways of thinking of embedding a probability
density function on $\R^n$ in a Hilbert space. The first is to simply
assume that the probability density function is square integrable and
hence lies directly in $L^2(\R^n)$. The second is to use the fact
that a probability density function lies in $L^1(\R^n)$ and is
non-negative almost everywhere. Hence $\sqrt{p}$ will lie in
$L^2(\R^n)$.

For clarity we will write $L^2_D(\R^n)$ when we think of $L^2(\R^n)$
as containing densities directly. The $D$ stands for direct. We
write ${\cal D} \subset L^2_D(\R^n)$ where ${\cal D}$ is the set
of square integrable probability densities (functions with integral
$1$ which are positive almost everywhere).

Similarly we will write $L^2_H(\R^n)$ when we think of $L^2(\R^n)$ as
being a space of square roots of densities. The $H$ stands for
Hellinger (for reasons we will explain shortly). We will write ${\cal
  H}$ for the subset of $L^2_H$ consisting of square roots of
probability densities.

We now have two possible ways of formalizing the notion of a family of
probability distributions. In the next section we will define a smooth
family of distributions to be either a smooth submanifold of $L^2_D$
which also lies in ${\cal D}$ or a smooth submanifold of $L^2_H$ which
also lies in ${\cal H}$. Either way the families we discussed earlier
will give us finite dimensional families in this more formal sense.

The Hilbert space structures of $L^2_D$ and $L^2_H$ allow us to define
two notions of distance between probability distributions which we
will denote $d_D$ and $d_H$. Given two probability distributions $p_1$
and $p_2$ we have an injection $\iota$ into $L^2$ so one defines the
distance to be the norm of $\iota(p_1) - \iota(p_2)$. So given two
probability densities $p_1$ and $p_2$ on $\R^n$ we can define:

\begin{eqnarray*}
d_H(p_1,p_2) &=& \left( \int (\sqrt{p_1} - \sqrt{p_2})^2 \dif \mu \right)^{\frac{1}{2}} \\
d_D(p_1,p_2) &=& \left( \int (p_1 - p_2)^2 \dif \mu \right)^{\frac{1}{2}}.
\end{eqnarray*}

Here $\dif \mu$ is the Lebesgue measure. $d_H$ defines the
{\em Hellinger distance} between the two
distributions, which explains are use of $H$ as a subscript.
We will write $\langle \cdot, \cdot \rangle_H$ for the inner product
associated with $d_H$ and $\langle \cdot, \cdot \rangle_D$
or simply $\langle \cdot, \cdot \rangle$ for
the inner product associated with $d_D$.

In this paper we will consider the projection of the
conditional density of the true state
of the system given the observations
(which is assumed to lie in ${\cal D}$ or ${\cal H}$)
onto a submanifold. The notion of projection only makes sense with
respect to a particular inner product structure. Thus we can consider
projection using $d_H$ or projection using $d_D$. Each has advantages
and disadvantages.

The most notable advantage of the Hellinger metric is that the $d_H$
metric can be defined independently of the Lebesgue measure and its
definition can be extended to define the distance between measures
without density functions (see Jacod
and Shiryaev~\citep{jacod87a}).
In particular the Hellinger distance is indepdendent of the
choice of parameterization for $\R^n$. This is a very attractive
feature in terms of the differential geometry of our set up.

Despite the significant theoretical advantages of the $d_H$ metric,
the $d_D$ metric has an obvious advantage when studying mixture
families: it comes from an inner product on $L^2_D$ and so commutes
with addition on $L^2_D$. So it should be relatively easy to
calculate with the $d_D$ metric when adding distributions as happens in mixture
families. As we shall see in practice, when one performs concrete
calculations, the $d_H$ metric works well for exponential families and
the $d_D$ metric works well for mixture families.  

\subsection{The tangent space of a family of distributions}

To make our notion of smooth families precise we need to explain what
we mean by a smooth map into an infinite dimensional space.

Let $U$ and $V$ be Hilbert spaces and let $f:U \to V $ be a continuous
map ($f$ need only be defined on some open subset of $U$). We say that
$f$ is Fre\'chet differentiable at $x$ if there exists a bounded
linear map $A:U \to V$ satisfying:
\[ \lim_{h \to x} \frac{\| f(h) - f(x) - A h
  \|_V}{\|h\|_U}\]
If $A$ exists it is unique and we denote it by
${\mathrm D}f(x)$. This limit is called the Fre\'chet derivative of
$f$ at $x$. It is the best linear approximation to $f$ at $0$ in
the sense of minimizing the norm on $V$.

This allows us to define a smooth map $f:U \to V$ defined on an open
subset of $U$ to be an infinitely Fre\'chet differentiable map.  We
define an {\em immersion} of an open subset of $\R^n$ into $V$ to be a
map such that ${\mathrm D}f(x)$ is injective at every point where $f$
is defined. The latter condition ensures that the best linear
approximation to $f$ is a genuinely $n$ dimensional map.

Given an immersion $f$ defined on a neighbourhood of $x$, we can think
of the vector subspace of $V$ given by the image of ${\mathrm D}f(x)$
as representing the tangent space at $x$.

To make these ideas more concrete, let us suppose that $p(\theta)$ is
a probability distribution depending smoothly on some parameter
$\theta = (\theta_1,\theta_2,\ldots,\theta_m) \in U$ where $U$ is some
open subset of $\R^m$. The map $\theta \to p(\theta)$ defines a map
$i:U \to {\cal D}$. At a given point $\theta \in U$ and for a vector
$h=(h_1,h_2,\ldots,h_m) \in \R^m$ we can compute the
Fr\'echet derivative to obtain:
\[ {\mathrm D} i (\theta) h = \sum_{i=1}^m \frac{\partial p}{\partial
  \theta_i} h_i \]

So we can identify the tangent space at $\theta$ with the following
subspace of $L^2_D$:
\begin{equation}
\label{basisForD}
\span \{ \frac{\partial p}{\partial \theta_1},
\frac{\partial p}{\partial \theta_2}, \ldots, 
\frac{\partial p}{\partial \theta_m} \} \end{equation}

We can formally define a smooth $n$-dimensional family of probability
distributions in $L^2_D$ to be an immersion of an open subset of
$\R^n$ into ${\cal D}$. Equivalently it is a smoothly parameterized
probability distribution $p$ such that the above vectors in $L^2$ are
linearly independent.

We can define a smooth $m$-dimensional family of probability
distributions in $L^2_H$ in the same way. This time let $q(\theta)$ be
a square root of a probability distribution depending smoothly on
$\theta$. The tangent vectors in this case will be the partial
derivatives of $q$ with respect to $\theta$. Since one normally
prefers to work in terms of probability distributions rather than
their square roots we use the chain rule to write the tangent space
as:
\begin{equation}
\label{basisForH}
 \span \{ \frac{1}{2 \sqrt{p}} \frac{\partial p}{\partial \theta_1},
\frac{1}{2 \sqrt{p}} \frac{\partial p}{\partial \theta_2}, \ldots, 
\frac{1}{2 \sqrt{p}} \frac{\partial p}{\partial \theta_m} \}
\end{equation}

We have defined a family of distributions in terms of a single
immersion $f$ into a Hilbert space $V$. In other words we have defined
a family of distributions in terms of a specific parameterization of
the image of $f$. It is tempting to try and phrase the theory in terms
of the image of $f$. To this end, one defines an {\em embedded
  submanifold} of $V$ to be a subspace of $V$ which is covered by
immersions $f_i$ from open subsets of $\R^n$ where each $f_i$ is a
homeomorphisms onto its image. With this definition, we can state that
the tangent space of an embedded submanifold is independent of the
choice of parameterization.

\begin{figure}[htp]
\begin{centering}
\includegraphics{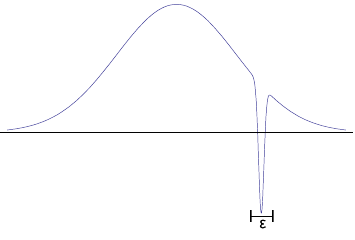}
\caption{An element of $L^2$ arbitrarily close to the normal distribution but not in $H$}
\label{fig:perturbednormal}
\end{centering}
\end{figure}

One might be tempted to talk about submanifolds of the space of
probability distributions, but one should be careful. The spaces
${\cal H}$ and ${\cal D}$ are not open subsets of $L^2_H$
and $L^2_D$ and so do not have any obvious Hilbert-manifold
structure. To see why, consider
Figure~\ref{fig:perturbednormal} where we have peturbed a probability
distribution slightly by subtracting a small delta-like function.

\subsection{The Fisher information metric}

Given two tangent vectors at a point to a family of probability
distributions we can form their inner product using $\langle \cdot, \cdot \rangle_H$.
This defines a so-called {\em Riemannian metric} on the family. With respect to a
particular parameterization $\theta$ we can compute the inner product
of the $i^{th}$ and $j^{th}$ basis vectors given in
equation~(\ref{basisForH}). We call this quantity $\frac{1}{4} g_{ij}$.
\begin{eqnarray*}
  \frac{1}{4} g_{ij}( \theta) & := & \langle \frac{1}{2 \sqrt{p}} \frac{ \partial
    p}{\partial \theta_i}, \frac{1}{2 \sqrt{p}} \frac{ \partial
    p}{\partial \theta_j} \rangle_H \\
  & = & \frac{1}{4} \int \frac{1}{p} \frac{ \partial
    p}{ \partial \theta_i} \frac{ \partial p}{ \partial \theta_j} \dif \mu \\
  & = & \frac{1}{4} \int \frac{ \partial
    \log p}{ \partial \theta_i} \frac{ \partial \log p}{ \partial \theta_j} p \dif \mu \\
  & = & \frac{1}{4} E_p( \frac{ \partial \log p}{\partial \theta_i} \frac{ \partial \log p}{\partial \theta_j})
\end{eqnarray*}
Up to the factor of $\frac{1}{4}$, this last formula is the standard
definition for the Fisher information matrix. So our $g_{ij}$ is the
Fisher information matrix. We can now interpret this matrix as the
Fisher information metric and observe that, up to the constant factor,
this is the same thing as the Hellinger distance. See \citep{amari85a},
\citep{murray93a} and \citep{aggrawal74a} for more in depth study on
this differential geometric approach to statistics.

\begin{example}The Gaussian family of densities can be parameterized
  using parameters mean $\mu$ and variance $v$. With this
  parameterization the Fisher metric is given by:
\[ g( \mu, v) = \frac{1}{v} \left[
\begin{array}{cc}
1 & 0 \\
0 & 1/(2v)
\end{array}
\right] \]
\end{example}

The representation of the metric as a matrix depends heavily upon the
choice of parameterization for the family.

\begin{example}The Gaussian family may be considered as a particular
  exponential family with parameters $\theta_1$ and $\theta_2$ given
  by:
\[ p(x,\theta) = \exp( \theta_1 x + \theta_2 x^2 - \psi(\theta) ) \]
where $\psi(\theta)$ is chosen to normalize $p$. It follows that:
\[ \psi(\theta) = \frac{1}{2} \log \left( \frac{\pi}{-\theta_2} \right) -
\frac{{\theta_1}^2}{ 4 \theta_2} \]
This is related to the familiar parameterization in terms of $\mu$ and
$v$ by:
\[ \mu = -\theta_1/(2 \theta_2), \quad v = \sigma^2 = (1/\theta_2 -
\theta_1^2/\theta_2^2)/2 \]
One can compute the Fisher information metric relative to the parameterization
$\theta_1$ to obtain:
\[ g(\theta) = \left[ \begin{array}{cc}
-1/(2 \theta_2) & \theta_1/(2 \theta_2^2) \\
\theta_1/(2 \theta_2^2) & 1/(2 \theta_2^2) - \theta_1^2/(2 \theta_2^3)
\end{array} \right] \]
\end{example}

The particular importance of the metric structure for this paper is
that it allows us to define orthogonal projection of $L^2_H$ onto the
tangent space.

Suppose that one has $m$ linearly independent vectors $w_i$ spanning
some subspace $W$ of a Hilbert space $V$. By linearity, one can write
the orthogonal projection onto $W$ as:
\[ \Pi(v) = \sum_{i=1}^m [\sum_{i=1}^m A^{ij} \langle v, w_j \rangle ] w_i \]
for some appropriately chosen constants $A^{ij}$. Since $\Pi$ acts as the identity
on $w_i$ we see that $A^{ij}$ must be the inverse of the matrix $A_{ij}=\langle w_i, w_j \rangle$.

We can apply this to the basis given in equation~(\ref{basisForH}). Defining
$g^{ij}$ to be the inverse of the matrix $g_{ij}$ we
obtain the following formula for projection, using the Hellinger metric,
onto the tangent space of a family of distributions:
\begin{eqnarray*}
\Pi_H(v) = \sum_{i=1}^m \left[ \sum_{j=1}^m 4 g^{ij} \langle
v, \frac{1}{2 \sqrt{p}} \frac{ \partial p}{ \partial \theta_j} \rangle_H \right] \frac{1}{2 \sqrt{p}} \frac{ \partial p}{ \partial \theta_i}
\end{eqnarray*}

\subsection{The direct $L^2$ metric}

The ideas from the previous section can also be applied to the
direct $L^2$ metric. This gives a different Riemannian metric on the manifold.

We will write $h=h_{ij}$ to denote the $L^2$ metric when written with respect to a particular parameterization.

\begin{example}
In coordinates $\mu$, $\nu$, the $L^2$ metric on the Gaussian family is:
\[ h( \mu, \nu) = \frac{1}{4 \nu \sqrt{ \nu \pi}} \left[
\begin{array}{cc}
1 & 0 \\
0 & \frac{3}{8 \nu}
\end{array}
\right]
\]
\end{example}

We can obtain a formula for projection in $L^2_D$ using the direct
$L^2$ metric using the basis given in equation~(\ref{basisForD}). We
write $h^{ij}$ for the matrix inverse of $h_{ij}$.
\begin{eqnarray}
\label{l2projectionformula}
\Pi_D(v) = \sum_{i=1}^m \left[ \sum_{j=1}^m h^{ij} \left\langle
v, \frac{ \partial p}{ \partial \theta_j} \right\rangle_D \right] \frac{ \partial p}{ \partial \theta_i}.
\end{eqnarray}

\section{The projection filter }\label{sec:PF}  

Given a family of probability distributions parameterised by $\theta$,
we wish to approximate an infinte dimensional solution to the non-linear filtering SPDE
using elements of this family. Thus we take the Kushner--Stratonovich
equation~(\ref{KSE:str}), view it as defining a stochastic vector field
in ${\cal D}$ and then project that vector field onto the tangent
space of our family. The projected equations can then be viewed as
giving a stochastic differential equation for $\theta$. In this section
we will write down these projected equations explicitly.

Let $\theta \to p(\theta)$ be the parameterization for our family. A
curve $t \to \theta( t)$ in the parameter space corresponds to a
curve $t \to p( \cdot , \theta(t))$ in ${\cal D}$. For such a curve, the
left hand side of the Kushner--Stratonovich equation~(\ref{KSE:str}) can
be written:

\begin{eqnarray*}
 d_t p(\cdot,\theta(t)) & = &
\sum_{i=1}^m \frac{\partial p(\cdot,\theta(t))}{\partial \theta_i}
d_t \theta_i(t) \\ 
& = & \sum_{i=1}^m v_i \dif \theta_i
\end{eqnarray*}
where we write $v_i=\frac{\partial p}{\partial \theta_i}$.
$\{v_i\}$ is the basis for the tangent space of the manifold at
$\theta(t)$.

Given the projection formula given in
equation~(\ref{l2projectionformula}), we can project the terms on the
right hand side onto the tangent space of the manifold using the
direct $L^2$ metric as follows:
\begin{eqnarray*}
\Pi_D^{\theta}[ {\cal L}^* p ] & = &
\sum_{i=1}^m \left[ \sum_{j=1}^m
h^{ij} \langle {\cal L}^* p ,  v_j \rangle \right] v_i \\
& = & 
\sum_{i=1}^m \left[ \sum_{j=1}^m
h^{ij} \langle p , {\cal L} v_j \rangle \right] v_i \\
\Pi_D^{\theta} [ \gamma^k( p ) ] & = &\sum_{i=1}^m \left[\sum_{j=1}^m
h^{ij} \langle \gamma^k( p ), v_j \rangle \right]
v_i
\end{eqnarray*}

Thus if we take $L^2$ projection of each side of equation~(\ref{KSE:str})
we obtain:
\[
\sum_{i=1}^m v_i \dif \theta^i = \sum_{i=1}^m
\left[
\sum_{j=1}^m h^{ij}
\left\{
\langle p, {\cal L}v_j \rangle \dif t - \langle \gamma^0(p), v_j \rangle \dif t
+ \sum_{k=1}^d \langle \gamma^k(p), v_j \rangle \circ \dif Y^k
\right\} \right] v_i
\]
Since the $v_i$ form a basis of the tangent space, we can equate the
coefficients of $v_i$ to obtain:
\begin{equation}
\label{KSE:l2projected}
\dif \theta^i = \sum_{j=1}^m h^{ij}
\left\{
\langle p(\theta), {\cal L}v_j \rangle \dif t - \langle \gamma^0(p(\theta)), v_j \rangle \dif t
+ \sum_{k=1}^d \langle \gamma^k(p(\theta)), v_j \rangle \circ \dif Y^k
\right\}.
\end{equation}
This is the promised finite dimensional stochastic differential
equation for $\theta$ corresponding to $L^2$ projection.

If preferred, one could instead project the Kushner--Stratonovich 
equation using the
Hellinger metric instead. This yields the following stochastic differential
equation derived originally in \citep{BrigoPhD}:
\begin{equation}
\label{KSE:hellingerProjected}
\dif \theta_i = \sum_{j=1}^m g^{ij}
\left( \langle \frac{{\cal L}^* p}{p}, v_j \rangle \dif t
     - \langle \frac{1}{2} |b|^2, v_j \rangle \dif t
     + \sum_{k=1}^d \langle b^k, v_j \rangle \circ \dif Y^k \right)
\end{equation}
Note that the inner products in this equation are the direct $L^2$ inner
products: we are simply using the $L^2$ inner product notation as a
compact notation for integrals.

Finally, it is now possible to explain why we resorted to the KS Equation for the optimal filter rather than the Zakai equation in deriving the projection filter.

\subsection{Projecting Kushner-Stratonovich vs Zakai}

Consider the nonlinear terms in the KS Equation (\ref{KSE:str}), namely
\[  \half\, p_t\, E_{p_t}\{\vert b_t \vert^2\} \,dt, \ \  \ 
    \sum_{k=1}^d p_t\, [-E_{p_t}\{b_t^k\}
] \circ dY_t^k .\]
Consider first the Hellinger projection filter (\ref{KSE:hellingerProjected}). By inspection, we see that there is no impact of the nonlinear terms in the projected equation.  Therefore, we see that projecting the Zakai unnormalized density SPDE in Hellinger metric onto the statistical manifold $p(\cdot,\theta)$  results in the same projection filter equation as projecting the normalized density KS SPDE. In other words, in Hellinger/Fisher metric the projection automatically takes care of normalization, so we could work with Zakai without changing the result. Intuitively, this happens because, under Hellinger, the set of square roots of normalized densities forms a unit sphere, and the tangent space is orthogonal to the direction corresponding to scaling (Hanzon 1987 \cite{hanzon87} uses Zakai rather than KS for this reason, among others). To see this, let us set $q_t(x) = \alpha_t p(x;\theta_t)$, so that $\alpha_t = \int q_t(x) dx$, and calculate ($d_t$ is short-hand for a Stratonovich differential)
\[\langle d_t \sqrt{q_t(x)}, \frac{\partial \sqrt{ p(x;\theta_t)}}{\partial \theta_i} \rangle =\frac{1}{4 \sqrt{\alpha}} \langle 
 (d \alpha_t) , \frac{\partial { p(x;\theta_t)}}{\partial \theta_i} \rangle  + \sqrt{\alpha} \sum_k    g_{k,i}(\theta)  \circ d \theta_k \   \]
and 
\[ \langle (d \alpha_t) , \frac{\partial { p(x;\theta_t)}}{\partial \theta_i} \rangle = (d \alpha) \frac{\partial}{\partial \theta_i} \int p(x;\theta) dx = (d \alpha) \  \frac{\partial 1}{\partial \theta_i} = 0\]
so that the component corresponding to change in the scaling constant does not contribute to the projection.

The main focus of this paper, however, is the $d_D$ projection filter (\ref{KSE:l2projected}). For this filter we do have an impact of the nonlinear terms. In fact, it is easy to adapt the derivation of the $d_D$ filter to the Zakai equation, which leads to the following new filter, namely the $d_D$ Zakai projection filter:
\begin{equation}
\label{ZAK:l2projected}
\dif \theta^i = \sum_{j=1}^m h^{ij}
\left\{
\langle p(\theta), {\cal L}v_j \rangle \dif t - \langle    \half\, \vert b_t \vert^2 \ p(\theta), 
  v_j \rangle \dif t
+ \sum_{k=1}^d \langle    b_t^k \ p(\theta), v_j \rangle \circ \dif Y^k
\right\}.
\end{equation}
This filter is clearly different from (\ref{KSE:l2projected}). We have implemented (\ref{ZAK:l2projected}) for the sensor case, using a simple variation on the numerical algorithms below. We found that (\ref{ZAK:l2projected}) gives slightly worse results for the $L^2$ residual of the normalized density than (\ref{KSE:l2projected}).
This can be explained simply by the fact that if we want to approximate $p$ in a given norm then we should project an equation for $p$ whereas if we wish to approximate $q$ we should project an equation for $q$. The fact that $\sqrt{q}$ has variable $L^2$ norm in time is not relevant for the projection in the Hellinger metric, while it is for the $d_D$ metric, where the lack of normalization in the Zakai Eq. density plays a role in the projection.

\section{Equivalence with ADFs and Galerkin Methods}\label{sec:Galerk}

The projection filter with specific metrics and manifolds can be shown to be equivalent to earlier filtering algorithms.
In particular, while the $d_H$ metric leads to the Fisher Information
and to an equivalence between the projection filter and  Assumed Density Filters (ADFs) when using exponential families, see \citep{brigo99}, the 
$d_D$ metric for simple mixture families is equivalent to a Galerkin method, as we show now following the second named author preprint \citep{brigo12}.

For applications of Galerkin methods to Nonlinear filtering we refer for example to \citep{nowak}, \citep{germani}, \citep{itok1}, \citep{beard}.

The basic Galerkin approximation is obtained by approximating the exact solution of the filtering SPDE (\ref{KSE:str}) with a linear combination of basis functions $\phi_i(x)$, namely
\begin{equation}\label{eq:galerkin:approx} \tilde{p}_t(x) := \sum_{i=1}^{\ell} c_i(t) \phi_i(x) .
\end{equation}

Ideally, the $\phi_i$ can be extended to indices $\ell+1,\ell+2,\ldots,+\infty$ so as form a basis of $L^2$. 

The method can be sketched intuitively as follows. 
We could write the filtering equation (\ref{KSE:str}) as
\[   \langle  -  d {p}_t + {\cal L}_t^\ast\, {p}_t\,dt
   -  \gamma_t^0({p}_t)\,dt
   + \sum_{k=1}^d  \gamma_t^k({p}_t) \circ dY_t^k\ ,  \xi \rangle = 0 \]
for all smooth $L^2$ test functions $\xi$ such that the inner product exists.

We replace this equation  with the equation
\[   \langle  -  d \tilde{p}_t + {\cal L}_t^\ast\, \tilde{p}_t\,dt
   -  \gamma_t^0(\tilde{p}_t)\,dt
   + \sum_{k=1}^d  \gamma_t^k(\tilde{p}_t) \circ dY_t^k\ ,  \phi_j \rangle = 0, \ \ j=1,\ldots,\ell . \]

By substituting Equation (\ref{eq:galerkin:approx}) in this last equation, using the linearity of the inner product in each argument and by using integration by parts we obtain easily an equation for the combinators $c$, namely

\begin{eqnarray}\label{eq:galerkinmix} \sum_{i=1}^\ell \langle \phi_i,\phi_j \rangle d c_i =   
\sum_{i=1}^\ell \langle \phi_i, {\cal L} \phi_j \rangle c_i \ dt -  
\langle \gamma_t^0\left(\sum_{h=1}^\ell c_h \phi_h \right) , \phi_j \rangle \\ \nonumber
   + \sum_{k=1}^d \langle  \gamma_t^k\left(\sum_{h=1}^\ell c_h \phi_h \right), \phi_j \rangle \circ dY_t^k  
\end{eqnarray}

Consider now the projection filter with the following choice of the manifold. 
We use the convex hull of a set of basic $L^2$ probability densities $q_1,\ldots,q_{m+1}$, namely the basic mixture family
\begin{equation}\label{eq:basicmix}  \{ p(\theta) := \theta_1 q_1 + \ldots + \theta_m q_m + [1-(\theta_1 + \ldots + \theta_m)] q_{m+1}, \ \sum_{i=1}^m \theta_i < 1, \ \theta_i \ge 0 \ \mbox{for all} \ \ i  \} .
\end{equation}

We can see easily that tangent vectors for the $d_2$ structure and the matrix $h$ are
\[ v_j = q_j - q_{m+1}, \ \ \ h_{i,j} = \langle q_i-q_{m+1},q_j-q_{m+1}\rangle, \ \ i,j =1,\ldots,m. \]
so that the metric is constant in $\theta$. If we apply the $d_2$ projection filter Equation (\ref{KSE:l2projected}) with this choice of manifold, we see immediately by inspection that this equation  coincides with the Galerkin method Equation (\ref{eq:galerkinmix}) if we take 
\[\ell = m+1, \  c_i = \theta_i   \ \mbox{and} \  \phi_i = q_i - q_{m+1}   \ \mbox{for} \ i=1,\ldots,m, \ \mbox{and} \  c_{m+1} = 1, \ \phi_{m+1} = q_{m+1}.\] 
The choice of the simple mixture is related to a choice of the $L^2$ basis in the Galerkin method. A typical choice could be based on Gaussian radial basis functions, see for example \cite{kormann}. 

We have thus proven the following first main theoretical result of this paper:

\begin{theorem}
For simple mixture families (\ref{eq:basicmix}), the $d_2$ projection filter (\ref{KSE:l2projected}) coincides with a Galerkin method (\ref{eq:galerkinmix}) where the basis functions are the mixture components $q$.
\end{theorem}

However, this equivalence holds only for the case where the manifold on which we project is the simple mixture family (\ref{eq:basicmix}). More complex families, such as the ones we will use in the following,  will not allow for a Galerkin-based filter and only the $L^2$ projection filter can be defined there. Note also that even in the simple case (\ref{eq:basicmix}) our $L^2$ Galerkin/projection filter  will be different from the Galerkin projection filter seen for example in \cite{beard}, because we use Stratonovich calculus to project the Kushner-Stratonovich equation in $L^2$ metric. In \cite{beard} the Ito version of the Kushner-Stratonovich Equation is used instead for the Galerkin method, but since Ito calculus does not work on manifolds, due to the second order term moving the dynamics out of the tangent space (see for example \cite{brigo99b}), we use the Stratonovich version instead. The Ito-based and Stratonovich based Galerkin projection filters will therefore differ for simple mixture families, and again, only the second one can be defined for manifolds of densities beyond the simplest mixture family.

\section{Numerical Software Design }\label{sec:NImp}  

Equations~(\ref{KSE:l2projected}) and (\ref{KSE:hellingerProjected}) both
give finite dimensional stochastic differential equations that we hope
will approximate well the solution to the full Kushner--Stratonovich
equation. We wish to solve these finite dimensional equations
numerically and thereby obtain a numerical approximation to the
non-linear filtering problem.

Because we are solving a low dimensional system of equations we hope to
end up with a more efficient scheme than a brute-force finite
difference approach. A finite difference approach can also
be seen as a reduction of the problem to a finite dimensional system.
However, in a finite difference approach the finite dimensional
system still has a very large dimension, determined by the number of
grid points into which one divides $\R^n$. By contrast the finite
dimensional manifolds we shall consider will be defined by only a
handful of parameters.

%\section{Software design }\label{sec:Soft}

The specific solution algorithm will depend upon numerous choices:
whether to use $L^2$ or Hellinger projection; which family of
probability distributions to choose; how to parameterize that family;
the representation of the functions $f$, $\sigma$ and $b$; how to
perform the integrations which arise from the calculation of
expectations and inner products; the numerical method selected to
solve the finite dimensional equations.

To test the effectiveness of the projection idea, we have implemented
a \CPP engine which performs the numerical solution of the
finite dimensional equations and allows one to make various selections
from the options above. Currently our implementation is restricted to
the case of the direct $L^2$ projection for a $1$-dimensional state
$X$ and $1$-dimensional noise $W$. However, the engine does allow one
to experiment with various manifolds, parameteriziations and functions
$f$, $\sigma$ and $b$.

We use object oriented programming techniques in order to allow this
flexibility. Our implementation contains two key classes
\class{FunctionRing} and \class{Manifold}.

To perform the computation, one must choose a data structure to
represent elements of the function space. However, the most effective
choice of representation depends upon the family of probability
distributions one is considering and the functions $f$, $\sigma$ and
$b$. Thus the \CPP engine does not manipulate the data structure
directly but instead works with the functions via the \class{FunctionRing}
interface. A UML (Unified Modelling Language \citep{UML})
outline of the \class{FunctionRing} interface is given
in table~\ref{UML:FunctionRing}.

\begin{table}[htp]
\begin{centering}
{\sffamily
\begin{tabular}{|l|} \hline
FunctionRing \\ \hline 
  + add( $f_1$ : Function, $f_2$ : Function ) : Function \\
  + multiply( $f_1$ : Function, $f_2$ : Function ) : Function \\
  + multiply( $s$ : Real, $f$ : Function ) : Function \\
  + differentiate( $f$ : Function ) : Function \\
  + integrate( $f$ : Function ) : Real \\ 
  + evaluate( $f$ : Function ) : Real \\
  + constantFunction( $s$ : Real ) : Function \\
  \hline
\end{tabular}
}
\caption{UML for the \class{FunctionRing} interface}
\label{UML:FunctionRing}
\end{centering}
\end{table}

\begin{table}[htp]
\begin{centering}
{\sffamily
\begin{tabular}{|l|} \hline
Manifold \\ \hline
+ getRing() : FunctionRing \\
+ getDensity( $\theta$ ) : Function \\
+ computeTangentVectors( $\theta$  : Point ) : Function* \\
+ updatePoint( $\theta$  : Point, $\Delta \theta$  : Real* ) : Point \\
+ finalizePoint( $\theta$  : Point ) : Point \\
\hline
\end{tabular}
}
\caption{UML for the \class{Manifold} interface}
\label{UML:Manifold}
\end{centering}
\end{table}

The other key abstraction is the \class{Manifold}. We give a UML
representation of this abstraction in table~\ref{UML:Manifold}. For
readers unfamiliar with UML, we remark that the $*$ symbol can be read ``list''.
For example, the computeTangentVectors function returns a list of functions.

The 
\class{Manifold} uses some convenient internal representation for a point,
the most obvious representation being simply the $m$-tuple
$(\theta_1, \theta_2, \ldots \theta_m)$.
On request the \class{Manifold} is able to provide the density associated with any
point represented as an element of the \class{FunctionRing}.

In addition the \class{Manifold} can compute the tangent vectors at any
point. The \method{computeTangentVectors} method returns a list of
elements of the \class{FunctionRing}  corresponding to each of the
vectors $v_i = \frac{\partial p}{\partial \theta_i}$ in turn. If the point
is represented as a tuple
$\theta=(\theta_1, \theta_2, \ldots \theta_n)$, the method
\method{updatePoint} simply adds the components of the tuple $\Delta \theta$ to
each of the components of $\theta$. If a different internal representation
is used for the point, the method should make the equivalent change to
this internal representation.

The \method{finalizePoint} method is called by our algorithm at the
end of every time step. At this point the \class{Manifold}
implementation can choose to change its parameterization for the
state. Thus the \method{finalizePoint} allows us (in principle at least) to
use a more sophisticated atlas for the manifold than just a single
chart.

One should not draw too close a parallel
between these computing abstractions and similarly named mathematical
abstractions. For example, the space of objects that can be
represented by a given \class{FunctionRing} do not need to form a
differential ring despite the \method{differentiate} method. This is
because  the \method{differentiate} function will not be called
infinitely often by the algorithm below, so the functions in the ring
do not need to be infinitely differentiable.

Similarly the \method{finalizePoint} method allows the
\class{Manifold} implementation more flexibility than simply changing
chart. From one time step to the next it could decide to use a
completely different family of distributions. The interface even
allows the dimension to change from one time step to the next. We do
not currently take advantage of this possibility, but adapatively
choosing the family of distributions would be an interesting topic for
further research.

\subsection{Outline of the algorithm}

The \CPP engine is initialized with a \class{Manifold}
object, a copy of the initial \class{Point} and \class{Function}
objects representing $f$, $\sigma$ and $b$.

At each time point the engine asks the manifold to compute the tangent
vectors given the current point. Using the multiply and integrate
functions of the class \class{FunctionRing}, the engine can compute the
inner products of any two functions, hence it can compute the metric
matrix $h_{ij}$. Similarly, the engine can ask the manifold for the
density function given the current point and can then compute
${\cal L}^*p$. Proceeding in this way, all the coefficients of
$\dif t$ and $\circ \dif Y$ in equation~(\ref{KSE:l2projected}) can be
computed at any given point in time.

Were equation~(\ref{KSE:l2projected}) an It\^o SDE one could now
numerically estimate $\Delta \theta$, the change in $\theta$ over a
given time interval $\Delta$ in terms of $\Delta$ and $\Delta Y$, the
change in $Y$. One would then use the \method{updateState} method to
compute the new point and then one could repeat the calculation
for the next time interval. In other words, were
equation~(\ref{KSE:l2projected}) an It\^o SDE we could numerically solve the SDE using
the Euler scheme.

However, equation~(\ref{KSE:l2projected}) is a Stratonovich SDE so the
Euler scheme is no longer valid. Various numerical schemes for solving
stochastic differential equations are considered in \citep{
burrageburragetian} and \citep{kloedenplaten}. One of the simplest
is the Stratonovich--Heun method described in
\citep{burrageburragetian}. Suppose that one wishes to solve the SDE:
\[ \dif y_t = f(y_t) \dif t + g( y_t) \circ \dif W_t \]
The Stratonvich--Heun method generates an estimate for the solution
$y_n$ at the $n$-th time interval using the formulae:
\begin{eqnarray*}
Y_{n+1} & = & y_n + f(y_n) \Delta + g(y_n) \Delta W_n \\
y_{n+1} & = & y_n + \frac{1}{2} (f(y_n) + f(Y_{n+1}))
\Delta + \frac{1}{2}( g(y_n) + g(Y_{n+1})) \Delta W_n
\end{eqnarray*}
In these formulae $\Delta$ is the size of the time interval
and $\Delta W_n$ is the change in $W$. One can think of $Y_{n+1}$
as being a prediction and the value $y_{n+1}$ as being a correction. Thus
this scheme is a direct translation of the standard Euler--Heun scheme for
ordinary differential equations. 

We can use the Stratonovich--Heun method to numerically solve equation
~(\ref{KSE:l2projected}). Given the current value $\theta_n$ for the
state, compute an estimate for $\Delta \theta_n$ by replacing $\dif t$
with $\Delta$ and $\dif W$ with $\Delta W$ in
equation~(\ref{KSE:l2projected}). Using the \method{updateState} method 
compute a prediction $\Theta_{n+1}$. Now compute a second estimate for
$\Delta \theta_n$ using equation~(\ref{KSE:l2projected}) in the state
$\Theta_{n+1}$. Pass the average of the two estimates to the \method{
updateState} function to obtain the the new state $\theta_{n+1}$.

At the end of each time step, the method \method{finalizeState} is
called. This provides the manifold implementation the opportunity to
perform checks such as validation of the state, to correct the
normalization and, if desired, to change the representation it uses for
the state.

One small observation worth making is that the equation~(\ref{KSE:l2projected})
contains the term $h^{ij}$, the inverse of the matrix
$h_{ij}$. However, it is not necessary to actually calculate the
matrix inverse in full. It is better numerically to multiply both
sides of equation~(\ref{KSE:l2projected}) by the matrix $h_{ij}$ and then compute
$\dif \theta$ by solving the resulting linear equations directly. This
is the approach taken by our algorithm.

As we have already observed, there is a wealth of choices one could
make for the numerical scheme used to solve equation~(\ref{KSE:l2projected}),
we have simply selected the most convenient. The
existing \class{Manifold} and \class{FunctionRing} implementations could
be used directly by many of these schemes --- in particular those
based on Runge--Kutta  schemes. In principle one might also consider
schemes that require explicit formulae for higher derivatives such as
$\frac{\partial^2 p}{\partial \theta_i \partial \theta_j}$. In this
case one would need to extend the manifold abstraction to provide this
information.

Similarly one could use the same concepts in order to solve
equation~(\ref{KSE:hellingerProjected}) where one uses the
Hellinger projection. In this case the \class{FunctionRing}
would need to be extended to allow division. This would in turn complicate
the implementation of the integrate function, which is why we have not yet
implemented this approach.

\section{The case of normal mixture families}\label{sec:normmixcase}

We now apply the above framework to normal mixture families. 
Let ${\cal R}$ denote the space of functions which can be written as finite
linear combinations of terms of the form:
\[ \pm x^n e^{ a x^2 + b x + c} \]
where $n$ is non-negative integer and $a$, $b$ and $c$ are constants.
${\cal R}$ is closed under addition, multiplication and differentiation, so it
forms a differential ring.

We have written an implementation of \class{FunctionRing} corresponding
to ${\cal R}$. Although the implementation is mostly straightforward some
points are worth noting.

Firstly, we store elements of our ring in memory as a collection of tuples
$(\pm, a, b, c, n)$. Although one can write:
\[ \pm x^n e^{ a x^2 + b x + c } = q x^n e^{a x^2 + b x} \]
for appropriate $q$, the use or such a term in computer memory should
be avoided as it will rapidly lead to significant rounding errors. A small amount
of care is required throughout the implementation to avoid such rounding errors.

Secondly let us consider explicitly how to implement integration for this ring.
Let us define $u_n$ to be the integral of $x^n e^{-x^2}$.
Using integration by parts one has:

\[ u_n:= \int_{-\infty}^{\infty} x^n e^{-x^2} \dif x
= \frac{n-1}{2} \int_{-\infty}^{\infty} x^{n-2} e^{-x^2} \dif x = \frac{n-1}{2} u_{n-2} \]

Since $u_0 = \sqrt{\pi}$ and $u_1 = 0$ we can compute $u_n$
recursively. Hence we can analytically compute the integral of
$p(x) e^{-x^2}$ for any polynomial $p$. By substitution, we can now
integrate $p(x-\mu) e^{-(x-\mu)^2}$ for any $\mu$. By completing
the square we can analytically compute the integral of
$p(x) e^{a x^2 + b x + c}$ so long as $a<0$. Putting all this together
one has an algorithm for analytically integrating the elements of ${\cal R}$.

Let ${\cal N}^{i}$ denote the space of probability distributions that
can be written as $\sum_{k=1}^i c_k e^{a_k x^2 + b_k x}$ for some
real numbers $a_k$, $b_k$ and $c_k$ with $a_k<0$. Given a smooth curve
$\gamma(t)$ in ${\cal N}^{i}$ we can write:

\begin{equation*}
\gamma(t) = \sum_{k=1}^i c_k(t) e^{a_k(t) x^2 + b_k(t) x}. 
\end{equation*}
We can then compute:
\begin{eqnarray*}
\frac{ \dif \gamma }{\dif t}
 & = & \sum_{k=1}^i \left( \left( \frac{ \dif a_k}{ \dif t} x^2
                   + \frac{ \dif b_k}{ \dif t} x \right) c_k e^{a_k x^2 + b_k x}
                   + \frac{ \dif c_k}{ \dif t} e^{a_k x^2 + b_k x}                   
                   \right) \\
 & \in &  {\cal R}
\end{eqnarray*}

We deduce that the tangent vectors of any smooth submanifold of
${\cal N}^{i}$ must also lie in ${\cal R}$. In particular this means that
our implementation of \class{FunctionRing} will be sufficient to
represent the tangent vectors of any manifold consisting of finite
normal mixtures.

Combining these ideas we obtain the second main theoretical result of the paper.

\begin{theorem}
Let $\theta$ be a parameterization for a family of probability
distributions all of which can be written as a mixture of at most $i$
Gaussians. Let $f$, $a=\sigma^2$ and $b$ be functions in the
ring ${\cal R}$. In this case one can carry out the direct $L^2$
projection algorithm for the problem given by equation~(\ref{Lanc1-1})
using analytic formulae for all the required integrations.
\end{theorem}

Although the condition that $f$, $a$ and $b$ lie in ${\cal R}$ may
seem somewhat restrictive, when this condition is not met one could
use Taylor expansions to find approximate solutions, although in such case rigorous convergence results need to be established.

Although the choice of parameterization does not affect the choice of
\class{FunctionRing}, it does affect the numerical behaviour of the
algorithm. In particular if one chooses a parameterization with
domain a proper subset of $\R^m$, the algorithm will break down the moment the
point $\theta$ leaves the domain. With this in mind, in the numerical examples
given later in this paper we parameterize normal mixtures of $k$ Gaussians 
with a parameterization defined on the whole of $\R^n$. We describe this
parameterization below.

Label the parameters $\xi_i$ (with $1\leq i \leq k-1$), $x_1$, $y_i$ (with $2\leq i \leq k$)
and $s_i$ (with $1 \leq i \leq k$). This gives a total of $3k-1$ parameters. So we can write
\[ \theta=(\xi_1,\ldots, \xi_{k-1}, x_1, y_2, \ldots, y_{k}, s_1, \ldots, s_{k}) \]
Given a point $\theta$ define variables as follows:
\begin{eqnarray*}
\lambda_1 & = & \logit^{-1}(\xi_1) \\
\lambda_i & = & \logit^{-1}(\xi_i) (1 - \lambda_1 - \lambda_2 - \ldots - \lambda_{i-1}) \qquad (2 \leq i \leq k-1) \\
\lambda_k & = & 1 - \lambda_1 - \lambda_2 - \ldots - \lambda_{k-1} \\
x_i & = & x_{i-1} + e^{y_i} \qquad (2 \leq i \leq k) \\
\sigma_i &= & e^{s_i}
\end{eqnarray*}
where the $\logit$ function sends a probability $p \in [0,1]$ to its log odds,
$\ln(p/1-p)$.
We can now write the density associated with $\theta$ as:
\begin{equation*}
p(x) = \sum_{i=1}^k
\lambda_i \frac{1}{\sigma_i \sqrt{2 \pi}} \exp( -\frac{(x - x_i)^2}{2 \sigma_i^2} )
\end{equation*}

We do not claim this is the best possible choice of parameterization,
but it certainly performs better than some more na\"ive parameteriations
with bounded domains of definition. We will call the direct $L^2$
projection algorithm onto the normal mixture family given with this
projection the {\em L2NM projection filter}.

\subsection{Comparison with the Hellinger exponential (HE) projection algorithm}

A similar algorithm is described in \citep{brigo98,brigo99} for projection
using the Hellinger metric onto an exponential family. We refer to
this as the {\em HE projection filter}.

It is worth highlighting the key differences between our algorithm and
the exponential projection algorithm described in \citep{brigo98}.

\begin{itemize}

\item In \citep{brigo98} only the special case of the cubic sensor was
considered. It was clear that one could in principle adapt the
algorithm to cope with other problems, but there remained symbolic
manipulation that would have to be performed by hand. Our algorithm
automates this process by using the \class{FunctionRing} abstraction.

\item When one projects onto an exponential family, the stochastic term 
in equation~(\ref{KSE:hellingerProjected}) simplifies to a term with constant coefficients.
This means it can be viewed equally well as either an It\^o or Stratonovich SDE.
The practical consequence of this is that the HE algorithm
can use the Euler--Maruyama scheme rather than the Stratonvoich--Heun scheme
to solve the resulting stochastic ODE's. Moreover in this case the Euler-Maruyama
scheme coincides with the generally more precise Milstein scheme.

\item In the case of the cubic sensor, the HE
algorithm requires one to numerically evaluate integrals such as:

\[ \int_{-\infty}^{\infty} x^n \exp( \theta_1 + \theta_2 x + \theta_3 x^2
+ \theta_4 x^4) \dif x \]

where the $\theta_i$ are real numbers. Performing such integrals
numerically considerably slows the algorithm. In effect one ends up
using a rather fine discretization scheme to evaluate the integral and
this somewhat offsets the hoped for advantage over a finite
difference method.

\end{itemize}

\section{Numerical Results }\label{NumRes}  

In this section we compare the results of using the direct $L^2$
projection filter onto a mixture of normal distributions with other numerical methods.
In particular we compare it with:

\begin{enumerate}

\item A finite difference method using a fine grid which we term the {
\em exact filter}. Various convergence results are known
(\citep{kushner} and \citep{kushnerHuang})
for this method. In the simulations shown below we use a grid with $1000$
points on the $x$-axis and $5000$ time points. In our simulations we
could not visually distinguish the resulting graphs when the grid was
refined further justifying us in considering this to be extremely
close to the exact result. The precise algorithm used is as described
in the section on ``Partial Differential Equations Methods'' in chapter
8 of Bain and Crisan~\citep{crisan10}.

\item The {\em extended Kalman filter} (EK). This is a somewhat
heuristic approach to solving the non-linear filtering problem but
which works well so long as one assumes the system is almost linear.
It is implemented essentially by linearising all the functions in the
problem and then using the exact Kalman filter to solve this linear
problem - the details are given in \citep{crisan10}. The EK filter is
widely used in applications and so provides a standard benchmark. However,
it is well known that it can give wildly innaccurate results for non-linear problems
so it should be unsurprising to see that it performs badly for most of the
examples we consider.

\item The HE projection filter. In fact we have implemented a
generalization of the algorithm given in \citep{BrigoPhD} that can cope
with filtering problems where $b$ is an aribtrary polynomial, $\sigma$
is constant and $f=0$. Thus we have been able to examine the
performance of the exponential projection filter over a slightly wider
range of problems than have previously been considered.

\end{enumerate}

To compare these methods, we have simulated solutions of the equations
~(\ref{Lanc1-1}) for various choices of $f$, $\sigma$ and $b$. We have
also selected a prior probability distribution $p_0$ for $X$ and then
compared the numerical estimates for the probability distribution $p$
at subsequent times given by the different algorithms. In the examples below
we have selected a fixed value for the intial state $X_0$ rather than drawing
at random from the prior distribution. This should have no more impact upon
the results than does the choice of seed for the random number generator.

Since each of the approximate methods can only represent certain
distributions accurately, we have had to use different prior
distributions for each algorithm. To compare the two
projection filters we have started with a polynomial exponential distribution
for the prior and then found a nearby mixture of normal distributions.
This nearby distribution was found using a gradient search algorithm
to minimize the numerically estimated $L^2$ norm of the difference of
the normal and polynomial exponential distributions. As indicated earlier,
polynomial exponential distributions and normal mixtures are qualitatively similar
so the prior distributions we use are close for each algorithm.

For the extended Kalman filter, one has to approximate the prior distribution
with a single Gaussian. We have done this by moment matching. Inevitably this
does not always produce satisfactory results.

For the exact filter, we have used the same prior as for the $L^2$
projection filter.

\subsection{The linear filter}

The first test case we have examined  is the linear filtering problem.
In this case the probability density
will be a Gaussian at all times --- hence if we project onto the two
dimensional family consisting of all Gaussian distributions there
should be no loss of information. Thus both projection filters
should give exact answers for linear problems. This is indeed the
case, and gives some confidence in the correctness of the computer
implementations of the various algorithms.

\subsection{The quadratic sensor}

The second test case we have examined is the {\em quadratic sensor}.
This is problem~(\ref{Lanc1-1}) with $f=0$, $\sigma=c_1$ and
$b(x)=c_2 x^2$ for some positive constants $c_1$ and $c_2$. In this problem the
non-injectivity of $b$ tends to cause the distribution at any time to
be bimodal. To see why, observe that the sensor provides no
information about the sign of $x$, once the state of the system has
passed through $0$ we expect the probability density to become
approximately symmetrical about the origin.  Since we expect the
probability density to be bimodal for the quadratic sensor it makes
sense to approximate the distribution with a linear combination
of two Gaussian distributions.

\begin{figure}[htp]
\begin{centering}
\includegraphics{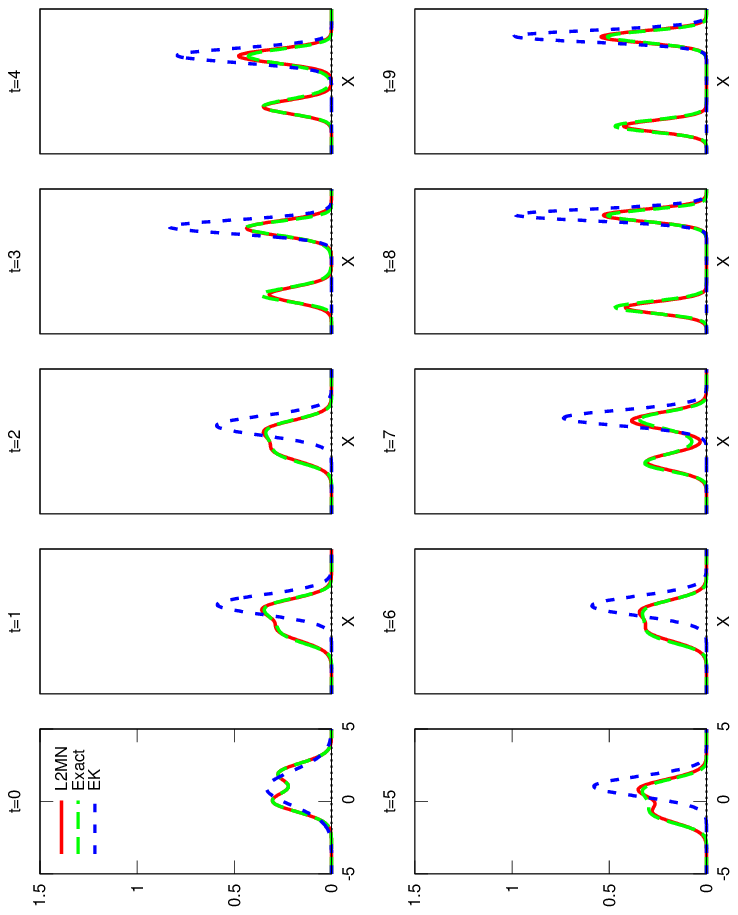}
\end{centering}
\caption{Estimated probability densities at $10$ time points for the problem $b(x)=x^2$}
\label{quadraticSensorTimePoints}
\end{figure}

In Figure~\ref{quadraticSensorTimePoints} we show the probability
density as computed by three of the algorithms at 10 different time
points for a typical quadratic sensor problem. To reduce clutter we have not
plotted the results for the exponential filter. The prior exponential
distribution used for this simulation was
$p(x)=\exp( 0.25 -x^2 + x^3 -0.25 x^4 )$. The initial state was $X_0 = 0$ and
$Y_0=0$.

As one can see  the probability
densities computed using the exact filter and the L2NM filter become
visually indistinguishable when the state moves away from the origin.
The extended Kalman filter is, as one would expect, completely unable
to cope with these bimodal distributions. In this case the extended Kalman filter
is simply representing the larger of the two modes.

\begin{figure}[htp]
\begin{centering}
\includegraphics{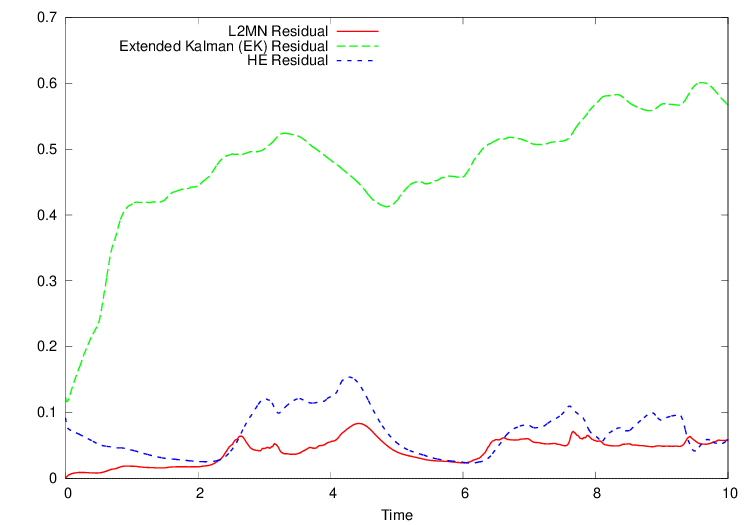}
\caption{$L^2$ residuals for the problem $b(x)=x^2$}
\label{quadraticSensorResiduals}
\end{centering}
\end{figure}

In Figure~\ref{quadraticSensorResiduals} we have plotted the {\em $L^2$
residuals} for the different algorithms when applied to the quadratic
sensor problem. We define the $L^2$ residual to be the $L^2$ norm of the
difference between the exact filter distribution and the estimated
distribution.
\begin{equation*}
\hbox{$L^2$ residual} = \left( \int | p_{\hbox{exact}} - p_{\hbox{approx}} |^2 \dif \mu \right)^{\frac{1}{2}}
\end{equation*}
As can be seen, the L2NM projection filter outperforms
the HE projection filter when applied to the quadratic sensor
problem. Notice that the $L^2$ residuals are initially small for both the
HE and the L2NM filter. The superior
performance of the L2NM projection filter in this case stems from the fact
that one can more accurately represent the distributions that occur
using the normal mixture family than using the polynomial exponential family.

If preferred one could define a similar notion of residual using the
Hellinger metric. The results would be qualitatively similar.

One interesting feature of Figure~\ref{quadraticSensorResiduals} is
that the error remains bounded in size when one might expect the error to
accumulate over time. This suggests that the arrival of new measurements
is gradually correcting for the errors introduced by the approximation.

\subsection{The cubic sensor}

A third test case we have considered is the {\em general cubic
sensor}. In this problem one has $f=0$, $\sigma=c_1$ for some constant
$c_1$ and $b$ is some cubic function.

The case when $b$ is a multiple
of $x^3$ is called the {\em cubic sensor} and was used as the test
case for the exponential projection filter using the Hellinger metric
considered in \citep{BrigoPhD}. It is of interest because it is the
simplest case where $b$ is injective but where it is known that the
problem cannot be reduced to a finite dimensional stochastic
differential equation \citep{HaMaSu}. It is known from earlier
work that the exponential filter gives excellent numerical results for
the cubic sensor.

\begin{figure}[htp]
\begin{centering}
\includegraphics{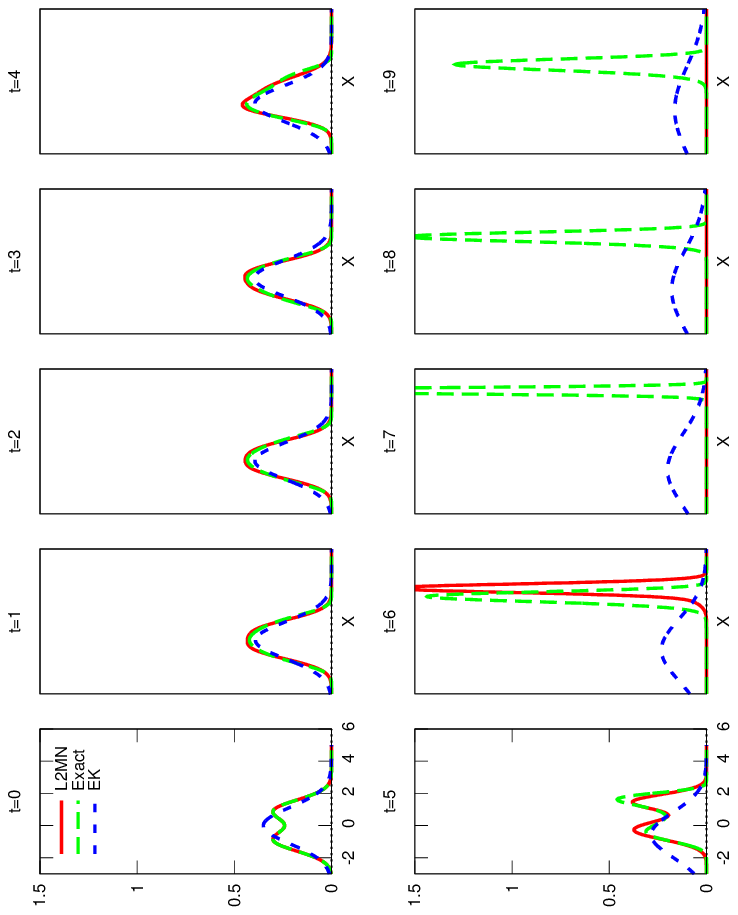}
\end{centering}
\caption{Estimated probability densities at $10$ time points for the problem $b(x)=x^3-x$}
\label{cubicSensorTimePoints}
\end{figure}

Our new implementations allow us to examine the general cubic sensor. In
Figure~\ref{cubicSensorTimePoints}, we have plotted example probability
densities over time for the problem with $f=0$, $\sigma=1$ and $b=x^3 -x$. 
With two turning points for $b$ this problem is very far from linear. As
can be seen in Figure~\ref{cubicSensorTimePoints} the L2NM projection
remains close to the exact
distribution throughout. A mixture of only two Gaussians is enough to
approximate quite a variety of differently shaped distributions with
perhaps surprising accuracy. As expected, the extended Kalman filter
gives poor results until the state moves to a region where $b$ is
injective. The results of the exponential filter have not been plotted
in Figure~\ref{cubicSensorTimePoints} to reduce clutter. It gave
similar results to the L2NM filter.

The prior polynomial exponential distribution used for this simulation
was $p(x)=\exp( 0.5 x^2 -0.25 x^4 )$. The initial state was $X_0 = 0$, which is one
of the modes of prior distribution. The inital value for $Y_0$ was taken to be $0$.

One new phenomenon that occurs when considering the cubic sensor
is that the algorithm sometimes abruptly fails. This is true for both
the L2NM projection filter and the HE projection filter.

\begin{figure}[htp]
\begin{centering}
\includegraphics{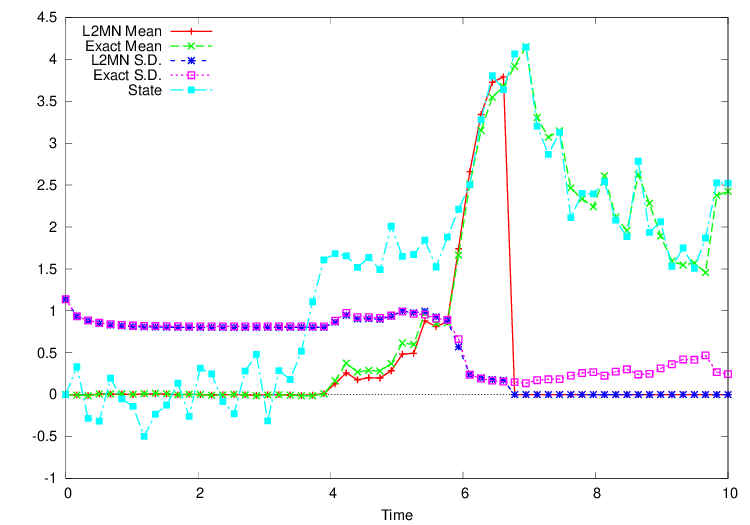}
\caption{Estimates for the mean and standard deviation for the problem $b(x)=x^3-x$}
\label{cubicSensorMeansAndSds}
\end{centering}
\end{figure}

To show the behaviour over time more clearly, in Figure~\ref{cubicSensorMeansAndSds}
we have shown a plot of the mean and standard deviation as estimated by the
L2NM projection filter against the actual mean and standard deviation.
We have also indicated the true state of the system. The mean for the L2MN filter
drops to $0$ at approximately time $7$. It is at this point that the
algorithm has failed.

What has happened is that as the state has moved to a region where the
sensor is reasonably close to being linear, the probability
distribution has tended to a single normal distribution. Such a
distribution lies on the boundary of the family consisting of a
mixture of two normal distributions. As we approach the boundary, $h_{ij}$
ceases to be invertible causing the failure of the algorithm. Analogous
phenomena occur for the exponential filter. 

The result of running numerous simulations suggests that the HE filter
is rather less robust than the L2NM projection filter. The typical behaviour
is that the exponential filter maintains a very low residual right up until
the point of failure. The L2NM projection filter on the other hand tends to give
slightly inaccurate results shortly before failure and can often correct itself
without failing.

\begin{figure}[htp]
\begin{centering}
\includegraphics{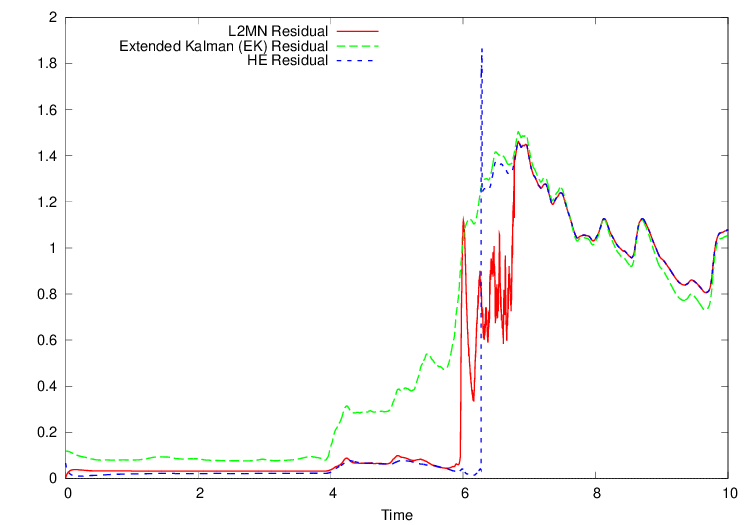}
\caption{$L^2$ residuals for the problem $b(x)=x^3-x$}
\label{cubicSensorResiduals}
\end{centering}
\end{figure}

This behaviour can be seen in Figure~\ref{cubicSensorResiduals}. In this
figure, the residual for the exponential projection remains extremely low until
the algorithm fails abruptly - this is indicated by the vertical dashed line. The L2NM
filter on the other hand deteriorates from time $6$ but only fails at time $7$.

The $L^2$ residuals of the L2MN method are rather large between times $6$ and $7$ but
note that the accuracy of the estimates for the mean and standard deviation in
Figure~\ref{cubicSensorMeansAndSds} remain reasonable throughout this time. To
understand this note that for two normal distributions with means a distance 
$x$ apart, the $L^2$ distance between the distributions increases as the standard
deviations of the distributions drop. Thus the increase in $L^2$ residuals
between times $6$ and $7$ is to a large extent due to the drop in standard
deviation between these times. As a result, one may feel that the $L^2$
residual doesn't capture precisely what it means for an approximation to be ``good''.
In the next section we will show how to measure
residuals in a way that corresponds more closely to the intuitive idea of them
having visually similar distribution functions. In practice one's definition of
a good approximation will depend upon the application.

Although one might argue that the filter is in fact behaving
reasonably well between times $6$ and $7$ it does ultimately fail.
There is an obvious fix for failures like this. When the
current point is sufficiently close to the boundary of the manifold,
simply approximate the distribution with an element of the boundary.
In other words, approximate the distribution using a mixture of fewer
Gaussians. Since this means moving to a lower dimensional family of distributions,
the numerical implementation will be more efficient on the boundary.
This will provide a temporary fix the failure of the algorithm,
but it raises another problem: as the state moves back into a region
where the problem is highly non linear, how can one decide how to leave the
boundary and start adding additional Gaussians back into the mixture?
We hope to address this question in a future paper.

\section{ Comparison with Particle Methods  }\label{sec:Part}  

Particle methods approximate the probability density $p$ using discrete
measures of the form:
	\[ \sum_i a_i(t) \delta_{v_i(t)} \]
These measures are generated using a Monte Carlo method. The measure can be
thought of as the empirical distributions associated with randomly located
particles at position $v_i(t)$ and of stochastic mass $a_i(t)$.

Particle methods are currently some of the most effective numerical
methods for solving the filtering problem. See \citep{crisan10} and the
references therein for details of specific particle methods and convergence
results.

The first issue in comparing projection methods with particle methods is
that, as a linear combination of Dirac masses, one can only expect
a particle method to converge weakly to the exact solution. In particular
the $L^2$ metric and the Hellinger metric are both inappropriate measures
of the residual between the exact solution and a particle approximation. Indeed
the $L^2$ distance is not defined and the Hellinger distance will always
take the value $\sqrt{2}$.

To combat this issue, we will measure residuals using the L{\'e}vy metric.
If $p$ and $q$ are two probability measures on ${\mathbb R}$ and $P$ and
$Q$ are the associated cumulative distribution functions then the
L{\'e}vy metric is defined by:
\[ d_L(p,q) = \inf \{ \epsilon:
P(x-\epsilon) - \epsilon \leq Q(x) \leq P(x + \epsilon) + \epsilon
\quad \forall x\} \]
This can be interpreted geometrically as the size of the largest square
with sides parallel to the coordinate axes that can be inserted between
the completed graphs of the cumulative distribution functions (the
completed graph of the distribution function is simply the graph of
the distribution function with vertical line segments added at
discontinuities).

The L{\'e}vy metric can be seen as a special case
of the L{\'e}vy--Prokhorov metric. This can be used to measure the
distance between measures on a general metric space. For Polish
spaces, the L{\'e}vy--Prokhorov metric metrises the weak convergence
of probability measures \citep{billingsley}. Thus the L{\'e}vy metric
provides a reasonable measure of the residual of a particle
approximation. We will call residuals measured in this way L{\'e}vy
residuals.

A second issue in comparing projection methods with particle methods
is deciding how many particles to use for the comparison. A natural
choice is to compare a projection method onto an $m$-dimensional
manifold with a particle method that approximates the distribution using
$\lceil (m+1)/2 \rceil$ particles. In other words, equate the dimension of the
families of distributions used for the approximation.

A third issue is deciding which particle method to choose for the
comparison from the many algorithms that can be found in the literature.
We can work around this issue by calculating the best possible
approximation to the exact distribution that can be made using
$\lceil (m+1)/2 \rceil$ Dirac masses. This approach will substantially underestimate the
L{\'e}vy residual of a particle method: being Monte Carlo methods, large
numbers of particles would be required in practice.

\begin{figure}[htp]
\begin{centering}
\includegraphics{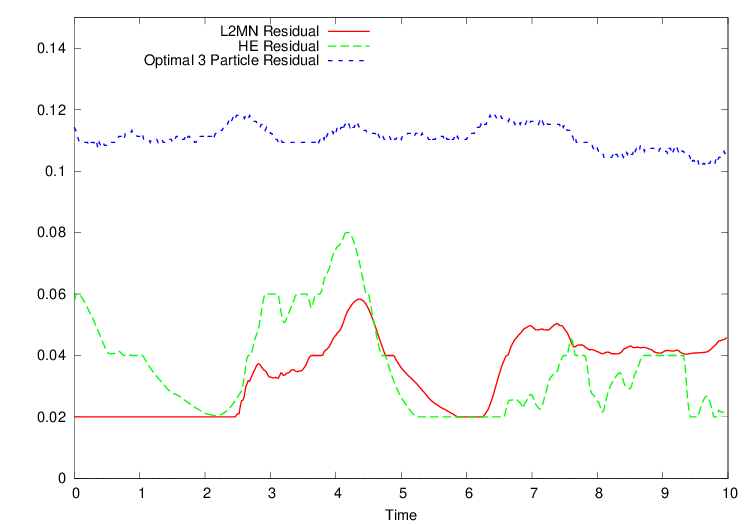}
\caption{L{\'e}vy residuals for the problem $b(x)=x^2$}
\label{levyResiduals}
\end{centering}
\end{figure}

In Figure \ref{levyResiduals} we have plotted bounds on
the L{\'e}vy residuals for the two projection methods for the
quadratic sensor. Since mixtures
of two normal distributions lie in a $5$ dimensional family, we have
compared these residuals with the best possible L{\'e}vy residual for
a mixture of three Dirac masses.

To compute the L{\'e}vy residual between two functions we have approximated
first approximated the cumulative
distribution functions using step functions. We have used the same
grid for these steps as we used to compute our ``exact'' filter.
We have then used a brute force approach to compute a bound on size of the largest
square that can be placed between these step functions. Thus if we have used
a grid with $n$ points to discretize the $x$-axis, we will need to make $n^2$
comparisons to estimate the L{\'e}vy residual. More efficient algorithms
are possible, but this approach is sufficient for our purposes.

The maximum accuracy of the computation of the L{\'e}vy metric is constrained
by the grid size used for our ``exact'' filter. Since the grid size
in the $x$ direction for our ``exact'' filter is $0.01$,
our estimates for the projection residuals are bounded below by $0.02$.

The computation of the minimum residual for a particle filter is a little
more complex. Let $\hbox{minEpsilon}(F,n)$ denote the minimum L{\'e}vy distance
between a  distribution with cumulative distribution $F$ and a distribution of
$n$ particles. Let $\hbox{minN}(F,\epsilon)$ denote the minimum number of particles
required to approximate $F$ with a residual of less than $\epsilon$. If
we can compute $\hbox{minN}$ we can use a line search to compute $\hbox{minEspilon}$.

To compute $\hbox{minN}(F,\epsilon)$ for an increasing step function $F$
with $F(-\infty)=0$ and $F(\infty)=1$, one needs to find the minimum
number of steps in a similar increasing step function $G$ that
is never further than $\epsilon$ away
from $F$ in the $L^{\infty}$ metric. One constructs candidate
step functions $G$ by starting with $G(-\infty)=0$ and then moving along
the $x$-axis adding in additional steps as required to remain within a distance
$\epsilon$. An optimal $G$ is found by adding in steps as late as possible and, when
adding a new step, making it as high as possible.

In this way we can compute $\hbox{minN}$ and $\hbox{minEpsilon}$ for step functions
$F$. We can then compute bounds on these values for a given distribution by
approximating its cumulative density function with a step function.

As can be seen, the exponential and mixture projection filters have similar
accuracy as measured by the L{\'e}vy residual and it is impossible to match
this accuracy using a model containing only $3$ particles.

\section{Conclusions}\label{sec:Conc} 

Projection onto a family of normal mixtures using the $L^2$ metric
allows one to approximate the solutions of the non-linear filtering
problem with surprising accuracy using only a small number of
component distributions. In this regard it behaves in a very similar
fashion to the projection onto an exponential family using the
Hellinger metric that has been considered previously.

The L2NM projection filter has one important advantage over the HE
projection filter, for problems with polynomial coefficients 
all required integrals can be calculated analytically. Problems
with more general coefficients can be addressed using Taylor series.
One expects this to translate into a better performing
algorithm --- particularly if the approach is extended to higher
dimensional problems.

We tested both filters against the optimal filter in simple but interesting systems, and
we provided a metric to compare the performance of each filter with the optimal one. 
We also tested both filters against a particle method, showing that with the same number of parameters the L2NM filter
outperforms the best possible particle method in Levy metric. 

We designed a software structure and populated it with models that make the L2NM filter
quite appealing from a numerical and computational point of view. 

Areas of future research that we hope to address include: the
relationship between the projection approach and existing numerical
approaches to the filtering problem; the convergence of the algorithm;
improving the stability and performance of the algorithm by adaptively
changing the parameterization of the manifold; numerical simulations
in higher dimensions.

%%%%%%%%%%%%%%%%%%%%%%%%%%%%%%%%%%%%%%
%
%
%   BIBLIOGRAPHY
%
%
%%%%%%%%%%%%%%%%%%%%%%%%%%%%%%%%%%%%%%

\end{document}